\documentclass[12pt]{article}
\usepackage{amsfonts,mathrsfs}
\def\d{{\rm d}}
\def\mi{{\rm i}}
\def\eps{\varepsilon}
\def\Re{\mathop{\rm Re\,}\nolimits}
\def\Im{\mathop{\rm Im\,}\nolimits}
\def\e{\mathop{\rm e}\nolimits}
\def\Res{\mathop{\rm Res}\nolimits}
\def\hf{{\textstyle{1 \over 2}}}
\def\qt{{\textstyle{1 \over 4}}}

\begin{document}

\title{A sharpening of Li's criterion\\
for the Riemann Hypothesis}
\author{{\bf Andr\'e Voros} \\
\\
CEA, Service de Physique Th\'eorique de Saclay\\
CNRS URA 2306\\
F-91191 Gif-sur-Yvette CEDEX (France)\\
{ E-mail : {\tt voros@spht.saclay.cea.fr}}\\
\\
and\\
\\
Institut de Math\'ematiques de Jussieu--Chevaleret\\
CNRS UMR 7586\\
Universit\'e Paris 7\\
2 place Jussieu,
F-75251 Paris CEDEX 05 (France)}

\maketitle

{\abstract
Exact and asymptotic formulae are displayed for the coefficients 
$\lambda _n$ used in Li's criterion for the Riemann Hypothesis.
In particular, we argue that if (and only if) the Hypothesis is true,  
$\lambda _n \sim n(A \log n +B)$ for $n \to \infty$ 
(with explicit $A>0$ and $B$). 
The approach also holds for more general zeta or $L$-functions.
}
\bigskip

Alternative expressions are presented for the coefficients 
$\lambda _n$ introduced by Li \cite{LI,BL}
to recast the Riemann Hypothesis as $ \lambda_n >0 \ (\forall n)$:
a first representation (\ref{LZS}) as a finite oscillatory sum,
then a closely related integral representation (\ref{LZI}),
upon which the saddle-point method finally
yields definite $n \to +\infty$ asymptotic estimates
in one of two sharply distinct forms, (\ref{RS}) or (\ref{REs}),
depending on the falsity or truth of the Riemann Hypothesis.
(Only the main ideas are indicated here.)

\section{Background.} 

The coefficients $\lambda_n$ are defined as
$ \lambda_n = \sum\limits _\rho \, [1-(1-1/\rho)^n] $,
or equivalently via the generating function
\begin{equation}
\label{LId}
{\d \over \d z} \log \Xi \Bigl( {1 \over 1-z}\Bigr) \equiv
\sum_{n=1}^\infty \lambda _n z^{n-1} , \qquad \qquad 
\bigl( \Xi(s) = s(s-1) \Gamma (s/2) \pi^{-s/2} \zeta(s) \bigr).
\end{equation}

The Riemann zeros are listed in pairs, as
\begin{equation}
\label{ZER}
\{\rho = \hf \pm \mi \tau_k \}_{k=1,2,\ldots}, \quad
\{ \Re \, \tau_k \} \mbox{ positive and non-decreasing;}
\end{equation}
sums and products over zeros are performed 
with $\rho$ and $(1-\rho)$ paired together, as usual;
we parametrize each such pair by the single number
$x_k=\rho(1-\rho)=\qt+{\tau_k}^2$.
A ``secondary" zeta function is
\begin{equation}
\label{ZDef}
Z(\sigma) = \sum_{k=1}^\infty {x_k}^{-\sigma}, \qquad \Re \sigma > \hf .
\end{equation}
It extends to a meromorphic function in $\mathbb C$, and
all its poles lie at the negative half-integers except $\sigma=+\hf$ \cite{K},
which has the polar part \cite{Vz}
\begin{equation}
\label{DPol}
Z(\hf+\eps) = R_{-2} \,\eps^{-2} + R_{-1} \,\eps^{-1} + O(1)_{\eps \to 0},
\qquad R_{-2} = {1 \over 8\pi} , \quad R_{-1} = -{\log 2\pi \over 4\pi}  .
\end{equation}

\section{Exact forms.} 

To reexpress the $\lambda _n$, we use 
a symmetrical Hadamard product form of $\Xi(s)$
\cite[Sec.1.10]{Ed}\cite{Vz}, namely 
$\Xi(s) = \prod_{k=1}^{\infty} \, \bigl[ 1 - s(1-s)/x_k \bigr]$,
instead of the standard form \cite[chap. 12]{Da}. Hence,
\begin{equation}
\label{HAL}
\Xi \Bigl( {1 \over 1-z} \Bigr) = 
\prod_{k=1}^{\infty} \biggl[ 1 + {z \over (1-z)^2 x_k} \biggr], \qquad
\log \Xi \Bigl( {1 \over 1-z} \Bigr) = 
- \sum_{j=1}^{\infty} {(-1)^j \over j} 
\biggl[ {z \over (1-z)^2 } \biggr] ^j Z(j);
\end{equation}
then, expanding $(1-z)^{-2j}$ by the generalized binomial formula,
reordering in powers of $z$ and identifying the output with (\ref{LId}), 
we get as first result
\begin{equation}
\label{LZS}
\lambda_n = - n \sum_{j=1}^n {(-1)^j \over j}{n+j-1 \choose 2j-1} \, Z(j) .
\end{equation}
Remark: other such linear relations already exist,
involving the cumulants of the Stieltjes constants \cite{BL} and/or the sums 
${\mathscr Z}_j = \sum_\rho \rho^{-j}$ instead of $Z(j)$;
see \cite[Sec.3.3]{Vz}  and references therein.
An advantage of (\ref{LZS}) is that the lesser known factors $Z(j)$ are
positive, and smooth (they sample the holomorphic function $Z(\sigma)$).
Still, this evaluation of $\lambda _n$ involves cancellations that increase
with~$n$.

An integral representation equivalent to (\ref{LZS}) simply by residue calculus,
but much more manageable, is
\begin{equation}
\label{LZI}
\lambda_n = {(-1)^n n \,\mi \over \pi} \int_C I(\sigma) \,\d \sigma,
\qquad I(\sigma) = 
{ \Gamma (\sigma+n) \Gamma (\sigma-n) \over \Gamma (2 \sigma+1) } \, Z(\sigma) ,
\end{equation}
where $C$ is a positive contour encircling (only) the subset of poles 
$\sigma=+1, \cdots, +n$ of the integrand $I(\sigma)$.

\section{Asymptotic forms.} 

The integral formula readily suggests
an asymptotic ($n \to \infty$) evaluation by the classic 
{\sl saddle-point method\/} \cite[Sec.2.5]{Er}.
First, the contour $C$ is deformed in the direction of decreasing $|I(\sigma)|$
until it crosses the nearest saddle-points of $|I(\sigma)|$. 
Among these saddle-points, the one(s) where $|I(\sigma)|$ is the largest 
give the dominant large-$n$ contributions to the integral, 
through local formulae at each such saddle-point
(the integrand itself may be asymptotically approximated as well:
e.g., by a Stirling formula for $\Gamma (z)$).
\medskip

{\sl We caution that our subsequent assertions, involving an integrand 
which cannot be described in fully closed form, partly retain an experimental
character and would warrant further confirmation.}
\medskip

In the present problem and for large $n$, the landscape of the function 
$|I(\sigma)|$ near the contour $C$ is dominantly controlled 
by the $\Gamma $-ratio: the resulting deformation 
is a dilation of $C$ away from the segment $[1,n]$, then the
nearest saddle-points can be of two types here
(once $n$ is large enough).

1) For $\sigma$ on the segment $(\hf,1)$,
$ |I(\sigma)| \sim \pi \,
[\sin \pi \sigma \, \Gamma (2 \sigma+1)]^{-1} n^{2\sigma -1} Z(\sigma)$
always has one {\sl real\/} minimum 
$ \sigma_{\rm r}(n) $ (tending to $\hf$ as $n \to \infty$;
other real saddle-points lie below $\sigma=\hf$, and are subdominant).

2) There may exist {\sl complex\/} saddle-points $\sigma_{\rm c}(n)$, 
with imaginary parts proportional to $n$; 
we may focus just on the upper half-plane, as the lower half-plane
gives complex-conjugate (``c.c.") contributions.
Wherever $Z(\sigma)$ is dominated by the term ${x_k}^{-\sigma}$
from one Riemann zero $x_k$ (naively, $x_1$), the saddle-point equation
(in the Stirling approximation for the $\Gamma $-ratio) is 
$0= {\d \over \d \sigma} \log |I(\sigma)| \sim \log(\sigma^2-n^2)-
2 \log 2\sigma -\log x_k$,
yielding the {\sl formal\/} saddle-point location
\begin{equation}
\label{CSP}
\sigma_{\rm c} (n) =  n \, \mi \, / \, 2 \tau_k .
\end{equation}

Here, any zero $x_k$ {\sl on\/} the critical axis ($\tau_k$ real) yields 
a purely imaginary $\sigma_{\rm c}(n)$, not eligible:
it lies outside the domain of convergence of (\ref{ZDef}),
and its contribution would be subdominant in any case. 
{\sl The discussion then fundamentally splits depending on the presence
or absence of zeros off the critical axis.}
\smallskip

[RH false] If the zero $x_k$ lies {\sl off\/} the critical axis 
(and selecting $\arg \tau_k >0$), the associated complex saddle-point 
$\sigma_{\rm c}(n)$ is relevant: 
it lies inside the domain of convergence $\{ \Re \sigma>\hf \}$
as soon as ${n \, |\Im 1/\tau_k| >1}$,
and it will exponentially dominate the real saddle-point (as seen later):
by the standard calculation, its asymptotic contribution to $\lambda _n$ is
$\bigl[ (\tau_k+\mi/2 )/( \tau_k-\mi/2) \bigr] ^n$, 
which {\sl grows exponentially\/} in modulus and fluctuates in phase.
This result can also be rigorously confirmed directly from (\ref{LId}):
by a conformal mapping argument \cite{BL},
the function ${\d \over \d z} \log \Xi \bigl( {1 \over 1-z}\bigr)$
has precisely the points $z_k={(\tau_k - \mi/2) (\tau_k + \mi/2)^{-1}}$ 
and $z_{-k}=\overline z_k$ as simple poles of residue 1 in the unit disk; 
by a general Darboux theorem \cite[chap. VII \S 2]{Di}, 
this implies the asymptotic form $\sum_k (z_{\pm k})^{-n}$ 
for the Taylor coefficients $\lambda _n$ of that function, i.e.,
\begin{equation}
\label{RS}
\lambda _n \sim \sum_{\{ \arg \tau_k >0 \}} 
\Bigl( {\tau_k+\mi/2 \over \tau_k-\mi/2} \Bigr) ^n
\quad [{}+ {\rm c.c.}] \quad \pmod {o(\e^{cn}) \ \forall c>0}, \quad 
n \to \infty. \qquad \quad [{\rm RH\ false}]
\end{equation}
(An infinite set of these zeros poses no extra problem:
their contributions form an asymptotic sequence.)

On the other hand, the Darboux approach does not resolve the case [RH true],
in which the poles $z_{\pm k}$ all lie at the same (unit) distance, 
and cluster at $z=1$~!
\smallskip

[RH true] By contrast, the saddle-point analysis of (\ref{LZI}) 
still appears to work. Now all the $\tau_k$ are real,
$Z(\sigma)=O(Z(\Re \sigma) \, |\Im \sigma|^{-3/2})$ in $\{ \Re \sigma > \hf \}$,
and the contour $C$ can be freely moved towards the boundary of the half-plane
without crossing any of the $\sigma_{\rm c}(n)$ (all purely imaginary).
Hence the only dominant saddle-point in this case is 
$\sigma_{\rm r}(n) \in (\hf,1)$;
it is shaped by the double pole of $Z(\sigma)$ at $\hf$
(itself generated by the totality of Riemann zeros), in the form
$ \sigma_{\rm r}(n) \sim \hf + {1 \over \log n}$.
One may then proceed as usual in the quadratic approximation of $\log I(\sigma)$
around $ \sigma_{\rm r}(n) $, but this is not so fit for a confluent case
($\sigma=\hf$ is a singular point).
Here, it is at once simpler and more accurate to
keep on deforming a portion of the contour $C$ nearest to $ \sigma=\hf$
until it fully encircles this pole (now clockwise), 
and to note that the resulting additions
to the integral are asymptotically negligible.
Hence for [RH true], the result is
\begin{equation}
\label{REs}
\lambda _n \sim (-1)^n 2n \Res_{\sigma=1/2} 
\Bigl[ { \Gamma (\sigma+n) \Gamma (\sigma-n) \over \Gamma (2 \sigma+1) } 
Z(\sigma) \Bigr]
\sim 2 \pi n \,\bigl[ 2 R_{- 2} \log n - 2 R_{-2} (1-\gamma ) + R_{-1} 
\bigl] \ ;
\end{equation}
($\gamma = $ Euler's constant); 
specifically for the Riemann zeros, using the explicit values (\ref{DPol}),
\begin{equation}
\label{RER}
\lambda _n \sim \hf n (\log n - \log 2\pi - 1+\gamma ),
\qquad n \to \infty. \qquad \qquad [{\rm RH\ true}]
\end{equation}
The form (\ref{REs}) covers more general situations: 
the zeros of, e.g., Dedekind zeta functions
and some Dirichlet $L$-functions \cite{VO}, in their [RH true] case
($R_{-2}$ is always positive, in agreement with Li's criterion).

We see a good agreement of (\ref{RER}) with numerical data \cite{M} 
for $n<3300$ -- still a bit short to give full confidence in (\ref{RER}), 
however. 
(This agreement is even improved in the mean if we include the contribution 
like (\ref{REs}) but from the next pole,
\begin{equation}
\delta \lambda _n = (-1)^n 2n \Res_{\sigma=0} 
\Bigl[ { \Gamma (\sigma+n) \Gamma (\sigma-n) \over \Gamma (2 \sigma+1) } 
Z(\sigma) \Bigr] = 2Z(0) =  +{7 \over 4},
\end{equation}
although this term should not be asymptotically meaningful -- larger 
oscillatory terms are also present.)
\bigskip

Finally, the result becomes even simpler for
some {\sl pure linear combinations of the Stieltjes cumulants\/}.
The latter may be defined in degree $n$ as
\begin{equation}
\label{SC}
g_n^{\rm c} = (-1)^{n-1} {\d^n \over \d s^n} [\log(s \, \zeta(1+s)]_{s=0} =
- \lim_{M \to +\infty} \Biggl\{ \sum_{m=1}^M {\Lambda(m) (\log m)^{n-1} \over m}
-{ (\log M)^n \over n} \Biggl\}
\end{equation}
(cf. formula (4.1) in \cite{BL}: this relates to the Euler factorization
of $\zeta(s)$ over the primes). Our $g_n^{\rm c}$ corresponds to
$(-1)^n (n-1)! \, \eta_{n-1}$ in \cite{BL}: e.g.,
$g_1^{\rm c}  = \gamma = -\eta_0$.
Theorem 2 in \cite{BL} evaluates the differences $(\lambda_n-S_n)$
where $S_n= -\sum\limits_{j=1}^n {n \choose j} \, \eta_{j-1} 
= \sum\limits_{j=1}^n {(-1)^{j-1} \over (j-1)!} {n \choose j}\,  g_j^{\rm c}$. 
As a consequence for $n \to \infty$, (\ref{RER}) implies
\begin{equation}
S_n = o(n) ; \qquad \qquad \qquad \qquad \qquad \qquad \qquad [{\rm RH\ true}]
\end{equation}
whereas in the opposite case, we expect $S_n \sim \lambda_n$ as in (\ref{RS}):
it will oscillate between {\sl exponentially large\/} values, 
negative {\sl and positive\/}, but appreciable in absolute size only beyond
$n \approx \min\limits_{\{ \arg \tau_k >0 \}} \{ [\Im 1/\tau_k]^{-1} \}$.

\end{document}